\documentclass[12pt]{article}

\usepackage[margin=1in]{geometry}  % set the margins to 1in on all sides
\usepackage{graphicx}              % to include figures
\usepackage{amsmath}               % great math stuff
\usepackage{amsfonts}              % for blackboard bold, etc
\usepackage{amsthm}                % better theorem environments
\usepackage{cite}

\usepackage{url} 

\usepackage{amssymb}
\usepackage{epstopdf}
\newtheorem*{claim}{Claim}

\title{The Geometry of Nim}
\author{Kevin Gibbons}
\date{}

\begin{document}
\maketitle

\begin{abstract}
We relate the Sierpinski triangle and the game of Nim. We begin by defining both a new high-dimensional analog of the Sierpinski triangle and a natural geometric interpretation of the losing positions in Nim, and then, in a new result, show that these are equivalent in each finite dimension.
\end{abstract}

\setcounter{section}{-1}
\section{Introduction}

\begin{figure}[b]
\begin{minipage}[b]{0.5\linewidth}
\centering
\includegraphics[scale=.3]{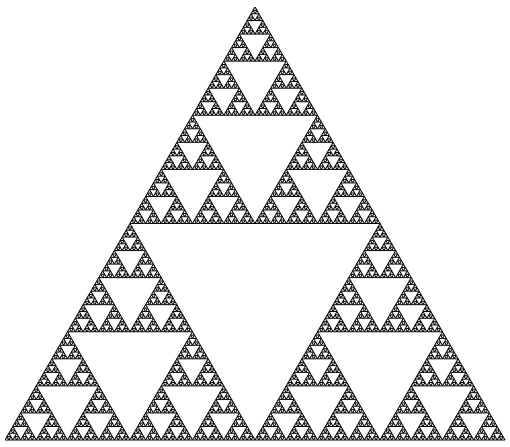}
\caption{The Sierpinski triangle}
\label{sierp}
\end{minipage}
\hspace{0.5cm}
\begin{minipage}[b]{0.5\linewidth}
\centering
\includegraphics[scale=.29]{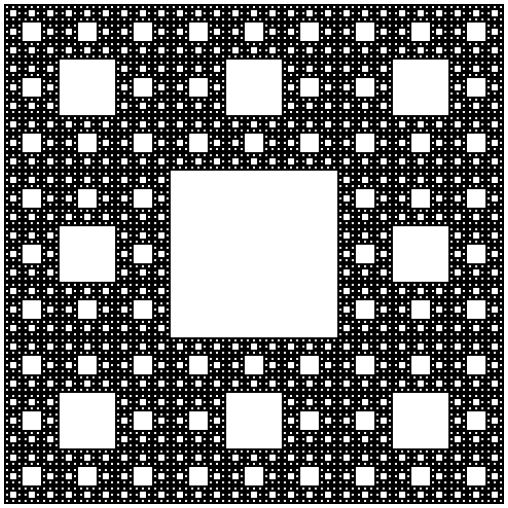}
\caption{A 2-dimensional Menger sponge}
\label{menger}
\end{minipage}
\end{figure}

The Sierpinski triangle (fig. \ref{sierp}) is one of the most recognizable figures in mathematics, and with good reason. It appears in everything from Pascal's Triangle to Conway's Game of Life. In fact, it has already been seen to be connected with the game of Nim, albeit in a very different manner than the one presented here\cite{Fraenkel}. A number of analogs have been discovered, such as the Menger sponge (fig. \ref{menger})
and a three-dimensional version called a tetrix. We present, in the first section, a generalization in higher dimensions differing from the more typical simplex generalization. Rather, we define a \emph{discrete Sierpinski demihypercube}, which in three dimensions coincides with the simplex generalization.

In the second section, we briefly review Nim and the theory behind optimal play. As in all impartial games (games in which the possible moves depend only on the state of the game, and not on which of the two players is moving), all possible positions can be divided into two classes - those in which the next player to move can force a win, called $N$-positions, and those in which regardless of what the next player does the other player can force a win, called $P$-positions. The second class is far less numerous, and it is those we consider here. When these positions appear in Nim is well understood, but we introduce a natural geometric interpretation, and formally define the set of such points in our interpretation.

Finally, we connect the two preceding sets by showing that the set of $P$-positions in Nim is precisely the Sierpinski demihypercube. This surprising fact is actually clear from the definitions we  introduce. We then briefly explore some conclusions and suggest questions for further investigation.

We recommend that the reader visit \url{http://www.kevingibbons.org/math/nim.html} for an in-browser rendering of the Sierpinski demihypercube, which can convey the structure better than the flat diagrams presented here.

\section{The Sierpinski Demihypercube}
\label{demihypercube}

The mouthful that is this section's title is in fact a fairly natural generalization of the Sierpinski triangle. A \emph{demihypercube} in geometry is a polytope (high-dimensional polyhedron) formed by removing alternating vertices of a cube and connecting those that remain. For example, the 2-demihypercube is merely a line segment: start with a square and remove two opposite vertices, and connect the remaining two. More interestingly, the 3-demihypercube is a tetrahedron (fig. 3).

\begin{figure}[b]
\begin{minipage}[b]{0.5\linewidth}
\centering
\includegraphics[scale=.3]{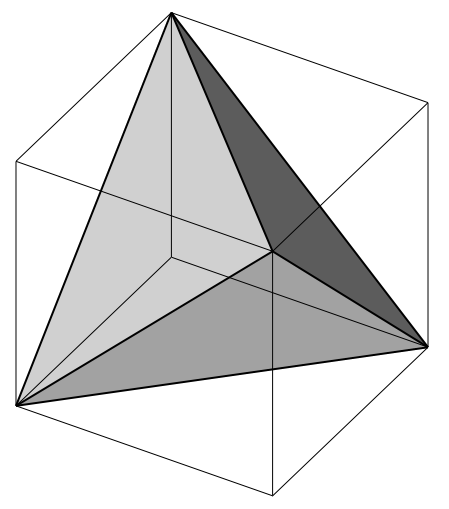}
\caption{The 3-demihypercube}
\end{minipage}
\hspace{0.5cm}
\begin{minipage}[b]{0.5\linewidth}
\centering
\includegraphics[scale=.15]{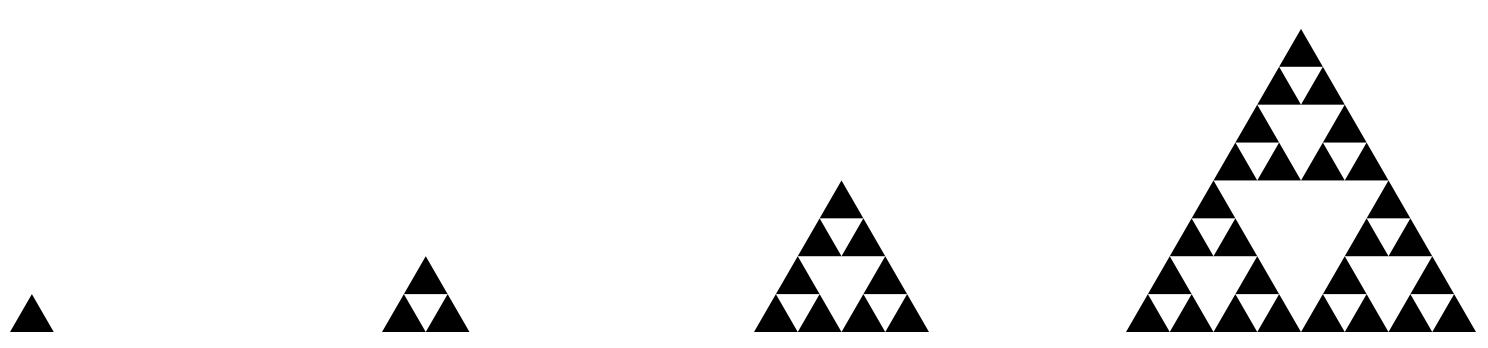}
\caption{Constructing a Sierpinski triangle}
\end{minipage}
\end{figure}

One can construct an approximation of the Sierpinski triangle in the following manner, pictured in fig. 4:
\begin{enumerate}
\item Start with an equilateral triangle, and define one vertex to be the origin.
\item Place a copy of the triangle with the origin-vertex at each of the original triangle's vertices.
\item Repeat this with the new figure.
\item Repeat for the desired number of iterations.
\end{enumerate}

We will construct the Sierpinski demihypercube in a very similar manner, except that we start with a discrete $d$-dimensional demihypercube, rather than a triangle. By \emph{discrete} we mean that we are only concerned with the vertices. To obtain the initial demihypercube, start at the origin of the unit hypercube and remove alternating vertices. This is equivalent to taking the set of points with all coordinates either zero or one, and with an even number of one coordinates. For example, the 3-dimensional discrete demihypercube is $\{(0,0,0), (1,1,0), (1,0,1), (0,1,1)\}$. Notice that these are the vertices of a tetrahedron. 

We are almost ready to formally define the general Sierpinski demihypercube, but first we need some notation. For a set $D \subseteq \mathbb{N}^n$ (where $\mathbb{N}$ is the set of nonnegative integers), we define
\begin{align*}
kD &= \{kx : x \in D\}\\
a+D &= \{ a + x : x \in D\} 
\end{align*}

Now we define an iteration of the discrete Sierpinski $d$-dimensional demihypercube recursively:
\begin{align*} 
D^1  &= \{ a \in \mathbb{N}^d : a \text{ has an even number of non-zero coordinates, all of which are 1}\}\\
D^{n+1} &= \bigcup_{a \in 2^{n}D^1} \left(a + D^{n}\right).
\end{align*}

This is practically the same as the process for generating the Sierpinski triangle described above: the base case is the initial demihypercube, and each iteration thereafter is given by placing a copy of the previous iteration at the desired locations. (The careful reader will note that we are placing the copies slightly offset from the vertices of the previous iteration. This is an artifact of making the demihypercube discrete.)

Finally, we define the full discrete Sierpinski demihypercube as the union of all iterations:
$$D = \bigcup_{n \in \mathbb{N}} D^n.$$

The first four iterations of the 3-dimensional version, $D^1$ through $D^4$, are presented below (fig. 5).
\begin{figure}[h]
  \begin{center}
    \includegraphics[scale=.25]{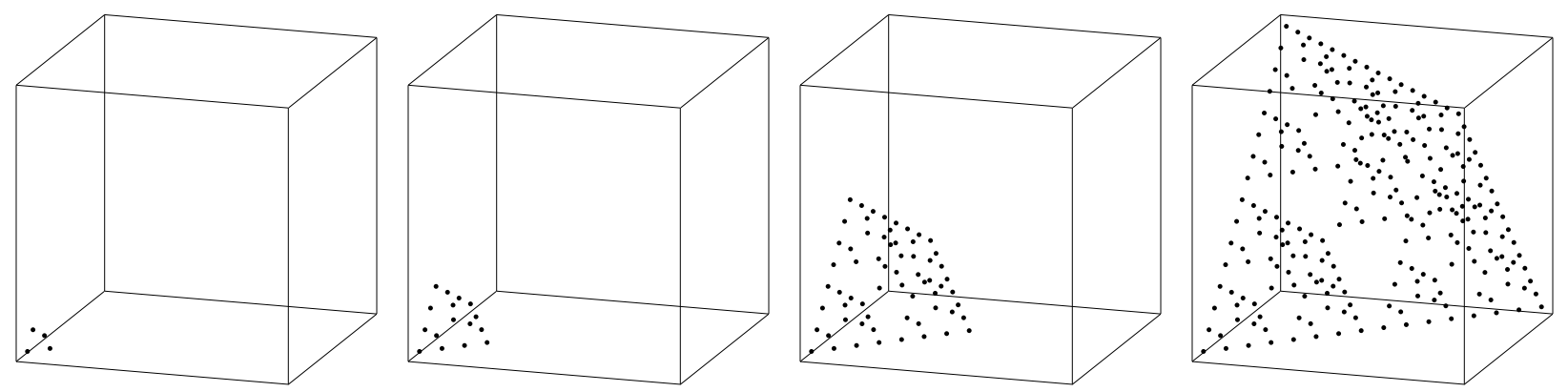}
  \end{center}
  \caption{The first four iterations of the discrete Sierpinski demihypercube}
\end{figure}

\section{Nim}

\subsection{An overview of the game and the theory}

To play Nim, find a friend and a bunch of stones. Divide the stones up into piles - as many piles, each of as many stones, as you want. You and your friend take turns removing from any pile as many stones (but at least one) as you would like. Whoever takes the final stone is the winner. This simple game is quite important in the study of combinatorial games - see \cite{winningways} for more.

There is a simple winning strategy for this game, one that may not be immediately obvious. To describe it, we first need to define the ``nim-sum'' operation on two natural numbers, written $x \oplus y$. It can be simply described as binary addition without carries. Formally, we write $x = x_02^0 + x_1 2^1 + x_2 2^2 + \dots + x_k2^k$, $y = y_0 2^0 + y_1 2^1 + y_2 2^2 \dots + y_k 2^k$ (where each $x_j$ and $y_j \in \{0, 1\}$). Then
$$x \oplus y = a_0 2^0 + a_1 2^1 + a_2 2^2 \dots + a_k 2^k,$$
where $a_j = x_j + y_j \mod 2$ (equivalently, $a_j$ is 1 if and only if $x_j$ or $y_j$, but not both, is  1). Computer scientists may recognize this as the ``bitwise exclusive or'' operation. (Since this operation is commutative and associative, we don't concern ourselves with order and omit parentheses.) To compute this efficiently, we can write all the numbers in binary and cancel pairs of equal powers, then add the remaining ones.

The winning strategy is to remove stones in such a way that the nim-sum over the number of stones in each pile is zero. (Say, for example, that there were piles of 4, 6, and 9 stones. The next player should reduce the pile of 9 to 2 stones, leaving 4, 6, and 2 stones respectively. Since the nim-sum of 4, 6, and 2 is 0, this is the correct move.) This is possible if and only if the nim-sum of the piles is not already zero; when it is the opposing player can force a win. It is beyond the scope of this paper to explain why this strategy is the correct one, but notice that the winning state, when all piles have zero stones, is such that the nim-sum of the piles is zero. The interested reader can convince herself that this strategy is indeed correct or consult \cite{Bouton} for a proof.

We will be concerned with those positions where the player not currently moving can force a win; these are called $P$-positions. By the previous paragraph, the set of $P$-positions is precisely the set of finite sets of natural numbers such that the nim-sum of the numbers in the set is zero. More specifically, when Nim is played with $d$ piles, the set of $P$-positions is the set of $d$-tuples of non-negative integers with nim-sum over all entries equal to zero.

\subsection{A geometric interpretation}

There is a very natural way to interpret a $d$-tuple of non-negative integers: as a point in $\mathbb{N}^d$. In this way we view the set of $P$-positions in $d$-pile Nim as a set of points in $d$-dimensional space. From the description above, we can write the set $P \subset \mathbb{N}^d$ of $P$-positions in $d$-pile Nim as
$$P = \{(a_1, a_2, \dots, a_d) \in \mathbb{N}^d : a_1 \oplus a_2 \oplus \cdots \oplus a_d = 0\}.$$

We need a final piece of machinery which will prove helpful later. We define the $n$-th restriction of $P \subseteq \mathbb{N}^d$ to be
$$P_n = \{(a_1, \dots, a_d) \in P : a_j < 2^n \text{ for all } j\}.$$

It is clear, if not immediately interesting, that
$$P = \bigcup_{n \in \mathbb{N}} P_n.$$

To make this interpretation concrete, below (fig. 6) are $P_1$ through $P_4$ in 3-pile Nim (notice that they are sets in three-space).

\begin{figure}[h] \label{asdf}
  \begin{center}
    \includegraphics[scale=.25]{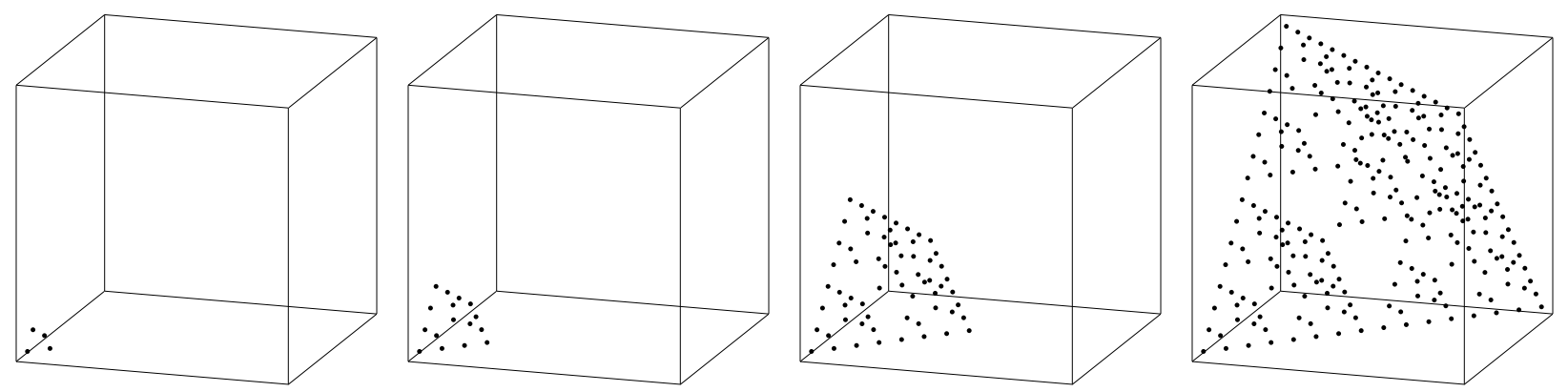}
  \end{center}
  \caption{The first four restrictions of the geometric interpretation of $P$-positions in Nim}
\end{figure}

This should seem familiar: it is in fact the same diagram used at the end of \textsection \ref{demihypercube}. In the following section, we show why.

\section{Putting It Together}

\subsection{What it's all about}

In the abstract, we claimed that the $P$-positions in Nim had a geometry similar to that of a generalized Sierpinski triangle. We are finally ready to describe and explain that claim.

\begin{claim} The set of $P$-positions in $d$-pile Nim is the full discrete Sierpinski $d$-demihypercube. \end{claim}

In the language introduced in the previous sections, $P = D$. As a reminder, we defined iterations of the Sierpinski demihypercube $D^n$ in the following manner:
\begin{align*} 
D^1  &= \{ a \in \mathbb{N}^d : a \text{ has an even number of non-zero coordinates, all of which are 1}\}\\
D^{n+1 }&= \bigcup_{a \in 2^{n}D^1} \left(a + D^{n}\right).
\end{align*}
Then $D$ was defined to be the union of all such $D^n$. The set $P$ of $P$-positions in $d$-pile Nim was seen to be 
$$P = \{(a_1, a_2, \dots, a_d) \in \mathbb{N}^d : a_1 \oplus a_2 \oplus \cdots \oplus a_d = 0\}.$$
After defining the $n$-th restriction of $P$ to be
$$P_n = \{(a_1, \dots, a_d) \in P : a_j < 2^n \text{ for all } j\},$$
we observed that $P$ was the union of all $P_n$.\\ \\

\begin{proof}
Intuitively, this is true because in the construction of $D^n$, shifting by some element in $2^{n}D^1$ is equivalent to adding coordinates with nim-sum zero, which does not affect the nim-sum of the resulting set.

Since
$$D = \bigcup_{n \in \mathbb{N}} D^n$$
and
$$P = \bigcup_{n \in \mathbb{N}} P_n,$$
it suffices to show that
$$D^n = P_n$$
for all $n \in \mathbb{N}$.

Formally, we prove this by induction.

When $n = 1$, we saw that $D^1$ was precisely those vertices of the unit hypercube with an even number of nonzero coordinates. This is $P_1$, since the nim-sum of piles which are zero or one is zero if and only if an even number of the piles are one.

For $n + 1 > 1$, we suppose that $P_{n} = D^{n}$.

We consider any point $x = (x_1, \dots, x_d)$ within the $2^{n+1}$ hypercube (which contains both $P_{n+1}$ and $D^{n+1}$). We write $x$ as the sum of its ``high part'' $h$ and its ``low part'' $y$, so that $x = h + y$. The high part $h = (h_1, \dots, h_d)$ is defined by $h_j = 2^{n}$ if $x_j \geq 2^{n}$, and 0 otherwise. That is, when $x$ is written as coordinates in binary, $h$ consists of the $2^n$ bits of the coordinates and $y$ of the remaining bits. For example, if $x = (7, 3, 5), n = 2$, we write $x = (4, 0, 4) + (3, 3, 1)$.

From the definition of nim-sum, we see that $x \in P_{n+1}$ if and only if (a) $h_1 \oplus h_2 \oplus \cdots \oplus h_d = 0$ and (b) $y \in P_{n}$. Similarly, we see $x \in D^{n+1}$ if and only if (a$'$) $h \in 2^{n}D^1$ and (b$'$) $y \in D^{n}$. But (a) is equivalent to (a$'$), since both just say that $h$ has an even number of nonzero coordinates, and (b) is equivalent to (b$'$) by the inductive hypothesis. Hence $P_n = D^n$ for all $n$.

\end{proof}

\subsection{Concluding remarks}

As always in mathematics, tying these two concepts together allows us to examine both in potentially enlightening ways. A couple are described below. We provide also some questions for further investigation.

The Sierpinski demihypercube is a highly symmetric figure, as one expects for a generalization of the Sierpinski triangle. For example, we believe its facets (high-dimensional faces) are ($d$-1)-dimensional Sierpinski demihypercubes and ($d$-1)-dimensional Sierpinski simplices. (A Sierpinski simplex is the truly natural high-dimensional generalization of the Sierpinski triangle, since the simplex is the natural high-dimensional generalization of the triangle.) Its properties are far from completely described here. 

The following is a non-obvious property implied by the result in this paper: for any specific $d$-1 coordinate values, there is precisely one point in the discrete Sierpinski demihypercube with those coordinates (the final coordinate varies depending on those coordinates). In particular, the final coordinate is given by the nim-sum of the specified coordinates, since the point with those coordinates has nim-sum zero and hence is a $P$-position in Nim. In other words, the shadow of the Sierpinski demihypercube onto a hyperplane perpendicular to an axis contains every lattice point between the coordinate axes precisely once. Fig. 7 shows this in the three-dimensional case.
\begin{figure}[h]
  \begin{center}
    \includegraphics[scale=.25]{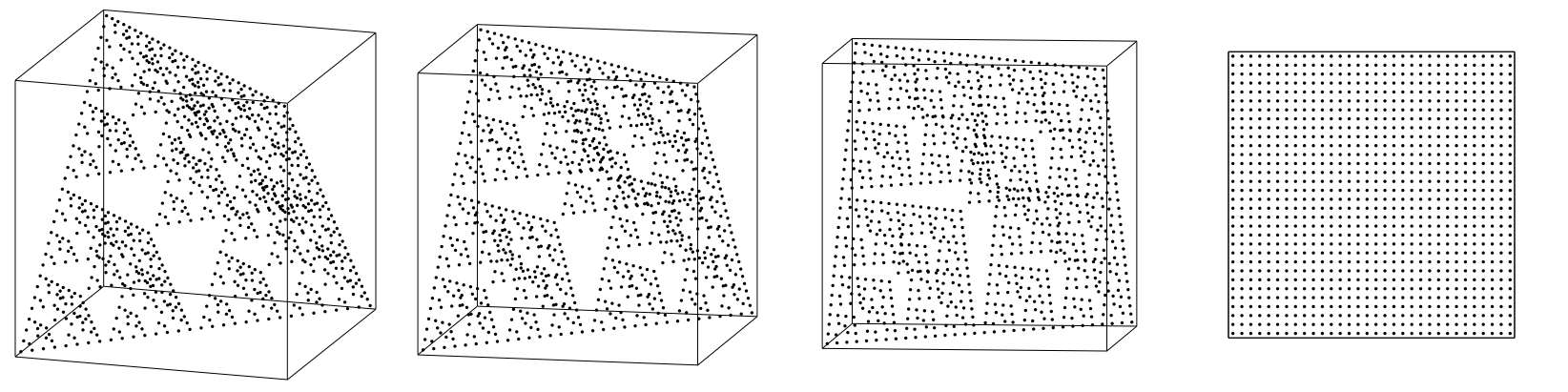}
  \end{center}
  \caption{Rotating the three-dimensional Sierpinski demihypercube to a side view}
\end{figure}

A connection between Nim and  the Sierpinski triangle which appears to us to bear no relation to the one presented here is described in \cite{Fraenkel}. It may be that both connections result from some deeper property which we have not seen.

Some generalizations of operations related to the nim-sum are presented in \cite{Fraenkel}. We suspect these generalizations also have some sort of high-dimensional Sierpinski structure.

Shadows of the Sierpinski demihypercube onto certain hyperplanes other than those perpendicular to coordinate axes have a Sierpinski-like symmetry. It may be worth exploring these shadows.

\section{Acknowledgments}

We would like to thank Mort Brown at the University of Michigan, Ann Arbor for presenting a game which suggested the geometric interpretation of Nim, for editorial assistance, and for miscellaneous help. We would also like to thank Ethan J. Goldberg,
Kenneth C. Millet, and Maribel Bueno, all at the University of California, Santa Barbara, for editorial and miscellaneous help.

\bibliographystyle{plain}
\bibliography{Nim}

\end{document}